
\magnification=\magstep1
\parindent=0in

\def\al{\alpha}
\def\b{\beta}

\def\d{\delta}
\def\e{\epsilon}
\def\g{\gamma}
\def\l{\lambda}
\def\m{\mu}
\def\n{\nu}
\def\o{\omega}

\def\t{\theta}

\def\R{{\bf R}}
\def\C{{\bf C}}
\def\P{{\bf P}}
\def\db{{\bar\pl}}

\def\i{\infty}
\def\I{\int}
\def\p{\prod}
\def\s{\sum}

\def\sub{\subseteq}
\def\ra{\rightarrow}

\def\G{\Gamma}
\def\D{\Delta}

\def\O{\Omega}

\def\v{\vskip .1in}

\def\[{{\bf [}}
\def\]{{\bf ]}}

\def\pl{\partial}

\magnification=\magstep1
\parindent=0in

\def\al{\alpha}
\def\b{\beta}

\def\d{\delta}
\def\e{\epsilon}
\def\g{\gamma}
\def\l{\lambda}
\def\m{\mu}
\def\n{\nu}
\def\o{\omega}

\def\si{\sigma}
\def\t{\theta}

\def\R{{\bf R}}
\def\C{{\bf C}}

\def\i{\infty}
\def\I{\int}
\def\p{\prod}
\def\s{\sum}

\def\sub{\subseteq}
\def\ra{\rightarrow}

\def\G{\Gamma}
\def\D{\Delta}
\def\cO{{\cal O}}
\def\cM{{\cal M}}
\def\cH{{\cal H}}

\def\cZ{{\cal Z}}

\def\cG{{\cal G}}

\def\cJ{{\cal J}}

\def\cX{{\cal X}}
\def\Aut{{\rm Aut}}

\def\ddb{\pl\bar\pl}
\def\cL{{\cal L}}

\def\v{\vskip .1in}

\def\[{{\bf [}}
\def\]{{\bf ]}}

\def\pl{\partial}

\def\bt{{\bf t}}
\def\us{{\underline{s}}}

\def\cN{{\cal N}}
\def\cD{{\cal D}}

\def\mapright#1{\smash{
   \mathop{\longrightarrow}\limits^{#1}}}

\settabs 3 \columns

\v\v\v
\centerline{\bf SCALAR CURVATURE, MOMENT MAPS,
AND}
\centerline{{\bf THE DELIGNE PAIRING }
\footnote*{Supported in part by
the National Science Foundation under grants
DMS-98-00783, DMS-02-45371, and DMS-01-00410.}}
\v\v

\centerline{D.H. Phong${}^{\dagger}$ and Jacob Sturm${}^{\ddagger}$}

\v\v

\centerline{${}^{\dagger}$ Department of Mathematics}
\centerline{Columbia University, New York, NY 10027}
\v
\centerline{${}^{\ddagger}$ Department of Mathematics}
\centerline{Rutgers University, Newark, NJ 07102}

\v\v
{\bf \S 1. Introduction.}
\v
Let $X$ be a compact complex manifold and $L \ra X$ a positive
holomorphic line bundle.   Assume that
$\Aut(X,L)/\C^\times$ is  discrete, where $\Aut(X,L)$ is the group of
holomorphic automorphisms of the pair $(X,L)$.
Donaldson [D] has recently proved that if $X$ admits a metric
$\o\in c_1(L)$ of constant scalar curvature, then
$(X,L^k)$ is Hilbert-Mumford stable for $k$  sufficiently
large. Since K\"ahler-Einstein metrics have constant scalar
curvature,
this confirms in one direction the well-known conjecture
of Yau [Y1-4] which asserts that the existence of a K\"ahler-Einstein
metric
is equivalent to stability in the sense of geometric invariant
theory. Additional evidence for Yau's conjecture had been
provided earlier by Tian [T2-4], who showed that
the existence of constant scalar curvature metrics implies $K$-stability
and
$CM$-stability.

\v

Donaldson's proof consists of showing that constant scalar curvature
implies the existence of a ``balanced basis" of $H^0(L^k)$,
that is, a basis which imbeds $X$ into projective space and
which is orthonormal with respect to the Fubini-Study metric.
The existence of a balanced basis is known to
imply the Hilbert-Mumford stability
of the manifold $X$, due to an earlier theorem
of Zhang [Z] and Luo [L]. The general theme of
approximating metrics via suitable projective
imbeddings had been advocated by Yau over
the years.

\v
To obtain a balanced basis of $H^0(L^k)$,
Donaldson interprets
the balanced
condition as the simultaneous vanishing of two moment maps
on a certain infinite dimensional manifold $\cH_0$:
One moment map, $\mu_{\cG}$, corresponds to the action of an
infinite dimensional
gauge group $\cG$,
and the other, $\mu_{su}$, to the action of a finite dimensional
unitary group $SU(N+1)$. Using Lu's formulas [Lu] for the
Tian-Yau-Zelditch
[T1][Ze][C] expansion of the density of states,
Donaldson begins by constructing a basis $\us_0$ for which
$\m_D$ is small, where
$\mu_D$ is the restriction of $\m_{su}$ to a single
$SL(N+1)$ orbit.  The solution to the  gradient flow equation for
$|\m_{D}|^2$ on the manifold $\cH_0//\cG$ provides a continuous
family of bases $\us_t$, starting at $\us_0$, whose
limit is the desired zero of the moment map. The key step
in the proof is to show that this limit exists, and
this is established by proving
a certain lower bound on
the derivative $d\m_{D}$ of $\mu_{D}$. Donaldson
proves this lower bound by an intricate infinite-dimensional analysis,
and also formulates two conjectures for its sharp forms.

\v
In this paper we improve on Donaldson's estimate for the lower
bound of $d\m_D$. A key idea is
the identification of the K\"ahler structure defined in [D] with the
curvature
of a certain line bundle on the symplectic quotient
$\cH_0//\cG$, namely the Deligne pairing
$\cM=\langle \pi^*O(1),\cdots,\pi^*O(1)\rangle_{\cX/\cZ}$
(c.f. (4.9) below). The Deligne pairing $\cM$ had been
introduced by Zhang [Z] in his study of heights of semistable
varieties. Its curvature is given by an explicit formula due to
Deligne [De]. Using this formula
and estimates for the $\db$ operator,
we  avoid the infinite-dimensional
gauge group and obtain the desired sharp bounds
for the derivative of the moment map. There are indications
(see Remark 2 in \S 5)
that our bounds are optimal since that is the case when $X=\C P^1$,
for vector fields inside $\Aut(X,O(k))^\perp\sub su(k+1) $.
\v

The paper is organized as follows:
In \S2 and \S3
we recall some basic facts about Deligne pairings and moment maps.
 In \S4  we
show that Donaldson's K\"ahler form is
the curvature of
the Deligne pairing $\cM$
(Theorem 1)
and
in \S 5 we use the curvature formula of [De] to get the desired
improvement over Donaldson's bound (Theorem 2).

\v

\v\v
\S 2. {\bf The Deligne pairing }
\v

We recall some of the basic definitions and properties in [De] and [Z]:
Let $\pi:\cX \ra S$ be a flat projective morphism of integral
schemes of relative dimension~$n$. Thus for every $s \in S$,
the fiber $\cX_s$ is a projective variety in $\P^N$ of dimension
$n$.
Let $\cL_0, \cL_1,...\cL_n$ be line bundles on $\cX$. The
Deligne pairing is a line bundle on $S$, denoted
$\langle\cL_0,\cL_1,...\cL_n\rangle(\cX/S)$,
and defined as follows:
Let $U \sub S$ be a small open set and let
$l_i$ be a rational section of $\cL_i$ over $\pi^{-1}U$.
Assume that the $l_i$ are chosen in ``general
position": This means $\cap_i\  div(l_i)=\emptyset $
and for each $s$ and $i$ with
$s\in U$ and $0\leq i \leq N$, the fiber
 $\cX_s$ is not
contained in $div(l_i)$.
 Then for every $k$, the map
$\big(\cap_{i\not= k} \ div(l_i)\big) \ra S$ is finite:
For
every $s$,
$\big(\cap_{i\not= k} \ div(l_i)\big)\cap \cX_s \sub \cX_s$ is
a zero cycle $\s n(s)P(s) $. This means the $n(s)$ are
integers and the $P(s)$ are a finite set of points in $\cX_s$.

\v
Now we  define $\langle\cL_0,\cL_1,...\cL_n\rangle(\cX/S)$. Over a
small $U \sub S$, this line bundle is trivial and
generated by the symbol $\langle l_0,...,l_n\rangle$ where the
$l_i$ are chosen to be in general position. If
$l_i'$ is another set of rational sections in
general position, then
$\langle l_0',...,l_n'\rangle = \psi(s)\langle l_0,...,l_n\rangle$ for
some
nowhere vanishing function $\psi$ on $U$ which we must specify.
We do this one section at a time: Assume that $l_i = l_i'$
for all $i \not=k$. Assume as well that the
rational function $f_k = l_k'/l_k$ is well defined
and non-zero on
$\big(\cap_{i\not= k} \ div(l_i)\big)$. Then
$ \psi(s) \ = \ \p f_k(P(s))^{n(s)}
$.
\v
Let $\pi: \cX \ra S$ and $\cL_0, ..., \cL_n$ as above.
Let $l$ be a rational section of $\cL_n$. Assume
all components of $div(l)$ are flat over $S$.
Then we have  the following
{\it  induction formula: }
$$  \langle\cL_0,...,\cL_n\rangle(\cX/S) \ = \
\langle \cL_0,...,\cL_{n-1}\rangle(div(l)/S)
\eqno(2.1)
$$
Assume now that $\cX,S$ are defined over $\C$, that $\cX$ is a
smooth variety and
that $\cL_i$ is endowed with a smooth hermitian metric.
We now define by induction a hermitian metric on
$\langle\cL_0,...,\cL_n\rangle(\cX/S)$:
Let $c_1'(\cL_i)=-{\sqrt{-1}\over 2\pi}\ddb \log ||l||^2$ be the
normalized curvature form of $\cL_i$,
where $l$ is an invertible
section of $\cL_i$.
When $n=0$,
$\langle\cL\rangle_s \ =  \otimes_{p \in \pi^{-1}(s)  } \ \cL_p $
so we define
$$
|| \langle l_0\rangle||_s \ =
 \ \p_{p \in \pi^{-1}(s)  }
||l_0(p)||
\eqno(2.2)
$$
In general, we define
$$
\log ||\langle l_0,..., l_n\rangle|| \ = \
\log||\langle l_0,...,l_{n-1}\rangle (div(l_n)/S)|| \ + \
\I_{\cX/S} \log ||l_n|| \Lambda_{i=0}^{n-1}c'_1(\cL_i)
\eqno(2.3)
$$
where the integral is the fiber integral over $S$.
If we combine the induction formula
with the definition of the metric,
we immediately get the following isometry:
$$
\langle \cL_0,...,\cL_n\rangle(\cX/S) \ = \
\langle\cL_0,...,\cL_{n-1}\rangle(div(l)/S)\otimes \cO
\left(-\I_{\cX/S} \log ||l_n||
\Lambda_{i=0}^{n-1}c'_1(\cL_i)\right)
\eqno(2.4)
$$
where $\cO(f)$ denotes the trivial line bundle with
metric $||1|| = exp(-f)$.
In particular,
$$
\langle\cL_0,...,\cL_{n-1},\cL_n\otimes\cO(\phi)\rangle(\cX/S) \ = \
\langle\cL_0,...,\cL_n\rangle(\cX/S)\otimes \cO(E)
\eqno(2.5)
$$
where
$$  E \ = \ \I_{\cX/S}\    \phi\cdot
\p_{k <n} c'_1(\cL_k)
\eqno(2.6)
$$
Using induction we get the following  {\it change of metric formula: }
$$
\langle\cL_0\otimes
\cO(\phi_0),...,\cL_n\otimes \cO(\phi_n)\rangle(\cX/S) \ = \
\langle\cL_0,...,\cL_n\rangle(\cX/S)\otimes \cO(E)
\eqno(2.7)
$$
where
$$  E \ = \ \I_X\  \s_{j=0}^n  \phi_j\cdot
\p_{k <j} c'_1(\cL_k\otimes \cO(\phi_k))\cdot\p_{k > j} c_1'(\cL_k)
\eqno(2.8)
$$
Finally, we recall a formula for the curvature  which
is given
by Proposition 8.5 of [De]:
$$
c_1'\big(\langle \cL_0,...,\cL_N\rangle(\cX/S)\big) = \ \I_{\cX/S}
\Lambda_{i=0}^n c_1'(\cL_i)
\eqno(2.9)
$$
A general result of this type has also been obtained
by Tian [T4].
\v
\v

{\bf \S 3. The Moment Map}
\v
The basic facts which we require from the theory of moment maps
are the following (see e.g. [DK]).
Let $K$ be a compact Lie group acting on a K\"ahler manifold
$(V,\o,I)$. Then the complexified group $K^c$
satisfies $Lie(K^c)=Lie(K)\otimes \C$, and $K^c$ acts on $(V,I)$,
preserving the complex structure, but not the K\"ahler form
or the metric.
\v
A moment map for the action of $K$ is a
smooth function
$\n: V\ra Lie(K)$ which satisfies the identity
$$
d\langle \n,\xi \rangle_{Lie(K)} \ = \ \iota_{X_\xi}\o
\eqno(3.1)
$$
where, for $\xi\in Lie(K)$, $X_\xi=\si(\xi)$ is the vector
field on $V$ generating the infinitesimal action of $\xi$,
and $\langle , \rangle_{Lie(K)}$ is an invariant Euclidean
metric on $Lie(K)$. We say that $\n$ is equivariant if
it intertwines the action of $K$ on $V$ with the adjoint
action of $K$ on $Lie(K)$. Moment maps need not always
exist, but if $K$ is semi-simple, then there is a
unique equivariant moment map $\n: V\ra Lie(K)$.
\v
Let $(V,\o,I)$ be a K\"ahler manifold. We say that
$(V,\o,I)$ has
a ``line bundle in the background"
if there is a triple $(L,h,A)$ where $L$ is a holomorphic
line bundle on $V$, $h$ is a hermitian metric on $L$ and
$A$ is a unitary connection on $L$ whose curvature, $F_A$,
satisfies: $F_A=-i\o$. Such a structure exists if and only
if the form ${1\over \pi}\o$ represents an integral cohomology
class.
\v
Now let $(V,\o,I)$ be a K\"ahler manifold and
assume $(L,h,A)$ is a line bundle in the background.
Then there exists an equivariant moment map
$\n: V \ra Lie(K)$ if and only if the action of $K$
on $(V,\o,I)$ can be lifted to an action of $K$ on
$(L,h,A)$. For example, if we are given $\n$
then the associated action of $K$ is given infinitesimally
by the formula:
$$
\hat\si(\xi) \ = \ \widetilde{\si(\xi) } + \n(\xi)\bt
\eqno(3.2)
$$
where $\bt$ is the infinitesimal action of $U(1)$ and,
for $Y$ a vector field on $V$,
$\tilde Y $ is the horizontal lift to $L$ given by the connection.
 We thus get
as well an action of $K^c$ on the holomorphic bundle $(L,I)$,
given infinitesimally by:
$$
\hat\si(\Xi)\ = \ [\si(\xi_1)+I\si(\xi_2) ]\,{\widetilde {} }\ + \
[\n(\xi_1)+i\n(\xi_2)]\bt
\eqno(3.3)
$$
where $\xi_1,\xi_2\in Lie(K)$ and
$\Xi=\xi_1+i\xi_2\in Lie(K^c)$.
\v

Now let $\tilde\G\sub L$ be a fixed orbit for $K^c$ acting
on $L$. Then $\tilde \G$ is a smooth manifold which lies over an orbit
$\G\sub V$ (also a smooth manifold). Define
$  h: \tilde \G \ra \R
$
by
 $h(\g)=\ -\log|\g|^2$.
Let $Q=K^c/K$ and
fix $\g_0\in \tilde \G$. Define
$  H: Q \ra \R
$
by
$H(g) \ = \ h(g\cdot\g_0)$.
The derivatives of $H_{\xi}(t)$ are given by the following basic formulas
(see [DK], \S 6.5.2):
\v

{\bf Proposition. }\hfill\break
{ \it
1. If $\g\in \tilde \G$, then $\g$ is a critical
point of $h$ if and only if $\n(\pi(\g))=0$ .}
\hfill\break
{\it
2. For $\xi\in Lie(K)$ let $H_\xi(t)= H(\exp(it\xi))$ and $x=x(t)={\rm
exp}(it\xi)\cdot\g_0$. Then}
$$
H_\xi'(t) \ = \ 2\langle \n(\exp(it\xi)\cdot x_0), \xi\rangle
\eqno(3.4)
$$
$$
H_\xi''(t)\ = \ 2\langle \si_x(\xi), \si_x(\xi)\rangle\ = \
2\o(\si_x(\xi), \overline{\si_x(\xi}))
\eqno(3.5)
$$
\v

\v
We give a proof for the convenience of the reader:
If $\Xi\in Lie(K^c)$ then
$\hat\si(\Xi)$ is a smooth vector field on $\tilde \G$.
We claim that the Lie derivative
$\cL_{\hat\si(\Xi)} h$ is given by
$$
(\cL_{\hat\si(\Xi)} h)(\g)\ = \ 2\langle \n(x), \xi_2\rangle
\eqno(3.6)
$$
where $x=\pi(\g)\in V$ (here $\pi:L\ra V$).
Indeed, clearly $\cL_{\tilde X}h=0$ for any vector
field $X$ on $V$ (since $|\g|^2$ is infinitesimally
constant in the horizontal direction). Thus, (3.3) implies
$$(\cL_{\hat\si(\Xi)} h)(\g)\ = -{d\over dt} \log
\bigg|\exp(it[\n(\xi_1)+i\n(\xi_2)]\g\bigg|^2\ = \
-{d\over dt}\log\exp(-2t\n(\xi_2))
$$
which yields (3.6). Taking $\Xi=i\xi$ in (3.6) we get
(3.4). Differentiating one
more time we get
$$  H_\xi''(t) \ = \ 2\langle d\n(\si_x(i\xi)), \xi\rangle
\ = \ 2\langle \si_x^*\si_x(\xi), \xi\rangle\ = \
2\langle \si_x(\xi), \si_x(\xi)\rangle\
$$
and this proves (3.5).
Finally, to prove statement 1 in the Proposition, observe
that (3.6) implies that
 $(\cL_{\hat\si(\Xi)} h)(\g)|_{t=0}\ =0$ for all $\Xi$ if
and only if $\langle \n(x_0),\xi_2\rangle = 0$ for
all $\xi_2\in Lie(K)$, that is, if and only if
$\n(x_0)=0$.
\v

\v\v
{\bf \S 4. Donaldson's K\"ahler structure and the Deligne pairing}
\v
\v
In this section we show how one can recover Donaldson's
moment map construction using the Deligne pairing. More
precisely, we show that the K\"ahler structure which Donaldson
defines has a line bundle $\cM$ in the background, where
$\cM$
is given by a certain Deligne pairing.
Thus Donaldson's K\"ahler structure is given by the curvature
of $\cM$, which can be computed using the general curvature
formula (2.9) of Deligne.
\v

\v
Let
$X$ be a compact complex manifold,
$p:L\ra X$ a positive holomorphic line bundle.
If $k$ is sufficiently large,
then any ordered  basis
$\us=(s_0,...,s_N)$  of $H^0(L^k)$
defines an imbedding
$\iota_\us : X \ra \P^N$ by
$\iota_\us(x)= (s_0(x),\cdots , s_N(x))$.
\v
An ordered basis $\us$ also defines a canonical
isomorphism of holomorphic bundles
$$
\psi_\us: L^k \ra \iota_\us^*O(1)
\eqno(4.1)
$$
as follows: Let $O(1)$ be the hyperplane line bundle on
$\P^N$, $\ell\in L^k$ and let
$x=p(\ell)$. Choose $j$ such that $s_j(x)\not=0$.
Then there is a unique complex number $a_j$ such that
$\ell=a_js_j(x)$. Then
$\psi_\us(\ell)$ is a linear map
sending the line spanned by
$(s_0(x),...,s_N(x))$ to $\C$. It is  defined
 by the formula
$$({s_0(x)\over s_j(x)},...,{s_N(x)\over s_j(x)})\mapsto a_j
\eqno(4.2)
$$
Note that this definition is independent
of the choice of $j$, and $\psi_\us$ has the property:
$ \psi_\us^*\iota_\us^* \pi_i = s_i$, where $\pi_i\in H^0(O(1))$ is
projection onto the $i^{th}$ component.
\v
Let $\us=(s_0,...,s_N)$ be an ordered basis of $H^0(L^k)$,
and let $h_{FS}$ be the Fubini-Study metric on
$O(1)$. Thus, if $(z_0,...,z_N)\in\C^{N+1}$ represents
a point in $\P^N$, and if $\pi_i\in H^0(O(1))$ is the projection
onto the $i^{th}$ component, then
$||\pi_i||^2_{h_{FS}} = { |z_i|^2\over \s_j |z_j|^2 }$.
Define
$$h_\us=\psi_\us^*\iota_\us^*h_{FS}
\eqno(4.3)
$$
Then $h_\us$ is a metric on $L^k$ which depends on $\us$.
One way to characterize $h_\us$ is as follows:
$$
||\cdot ||^2_{h_{\us}} \ = \
{||\cdot ||^2_{h_0}\over \s_j || s_j ||^2_{h_0}
}
\eqno(4.4)
$$
where $h_0$ is any fixed metric on $L^k$.
\v

\v
Now fix $k$,
a large positive integer and define
$$\tilde \cZ\ = \
 \{  \us=(s_0,...,s_N) : \{s_0,...,s_N\}
\ \hbox{is a basis of\  } H^0(L^k)\ \}/\C^\times
\eqno(4.5)
$$
and
$$  \cZ \ = \ \tilde \cZ/(PAut(X,L^k))
\eqno(4.6)
$$
where $PAut(X,L^k)=Aut(X,L^k)/\C^\times$
and the map $Aut(X,L^k) \hookrightarrow GL(N+1)$
is given by the action on global sections.
\v
Assume now that $PAut(X,L^k)$ is discrete. Then the construction
(4.5-4.6) provides $\cZ$ with a complex structure.
Next Donaldson defines a map $\m_D: \cZ \ra su(N+1)$ by
$$
\m_D(\us)\ = \ P_{su(N+1)}\big[i\langle s_\al, s_\b\rangle_{h_\us}\big]
\eqno(4.7)
$$
where $P_{su(N+1)}: u(N+1)\ra su(N+1)$ is the projection
onto the trace free subspace:
$P_{su(N+1)}(C_{\al\b})= C_{\al\b}-
[{1\over N+1 }tr(C)]\cdot \d_{\al\b}$
(the Lie algebra $su(N+1)$ is identified
with its dual using the invariant Hilbert-Schmidt pairing).
A key property of $\m_D$
is that it is the moment map for a certain symplectic structure
$\O_D$ on $\cZ$ which Donaldson constructs as follows.

\v
Fix a hermitian metric $h_0$
on $L$ and let $I_0$ be the holomorphic structure
on $X$. Let $A$ be the unitary connection on $L$ which
is compatible with $I_0$.  Let $\cJ_{int}$ be
the set of all integrable holomorphic structures on $X$
(so that in particular, $I_0\in \cJ_{int}$) and define
$\cH_0\sub \G(L^k)\times\cdots \times \ \G(L^k)\times \cJ_{int}$
to be the set of pairs $(\us, I)$ where $\us=(s_0,...,s_N)$ is
an ordered basis of $H^0((L,I)^k)$. Here $(L,I)$ is the
holomorphic line bundle determined by $I$ and the connection $A$.
Let $\cG$ be the group of $C^\i $
automorphisms of the triple $(L,h_0,A)$, so
that
$Lie(\cG)=C^\i(X)/\R$.
Then $\cG$ acts on $\cH_0$ and it turns out that there is a moment map
$\m_\cG: \cH_0\ra C^\i(X)$ for the action of $\cG$ (it is not
obvious that such a moment map should exist since $\cG$ is
infinite dimensional). Although the complexification $\cG^c$
of $\cG$ does not exist, one can, for $w\in \cH_0$, still consider the
``orbit" $\cG^c\cdot w \sub \cH_0$. Fix any real number $a>0$.
One shows
that every orbit meets $\m_\cG = a$ uniquely, up to the action
of $\cG$. In this way, one obtains, by symplectic reduction,
a K\"ahler structure on
$\cH_0//\cG = \cH_0/\cG^c$, and a corresponding
moment map which is the restriction to the set $\{\m_\cG=a\}$
of the map $\mu_{su}$ given by
$\m_{su}(\us)=P_{su(N+1)}\big[i\langle s_\al, s_\b\rangle_{h_0}\big]$.
The natural map $\b: \cZ\ra \cH_0//\cG$ is an imbedding, and provides
$\cZ$ with
the desired K\"ahler structure $\O_D$.  Moreover, the moment map
$\mu_D$ satisfies the relation: $\m_D=\m_{su}\circ \b$.
\v

\v\v

We now give a simple description of $\O_D$: Let
$$\tilde \cX \ = \ \{(x,\us): x\in \P^N,\  \us=(s_0,...,s_N)
\ \hbox{a basis of\  } H^0(L^k),\  x \in \iota_\us(X) \}
\eqno(4.8)
$$
and let $\cX = \tilde\cX/PAut(X,L^k)$. Then
$\cX\ra \cZ$ is a smooth holomorphic fibration whose
fibers are all isomorphic to $X$.
Let $\pi:\cX\ra \P^N$
 be  projection onto the second factor, and let
$$
\cM \ = \ \langle \pi^*O(1),...,\pi^*O(1) \rangle{(\cX/\cZ )}
\eqno(4.9)
$$
be the Deligne pairing of $n+1$ copies of the line bundle
$\pi^*O(1)$, in the sense explained in \S 2.
The line bundle $\cM$
is a hermitian line bundle on $\cZ$, invariant under
the action of $SU(N+1)$.  Let $\O_\cM$ be the curvature
of $\cM$. Formula (2.9) says that $\O_\cM$, the curvature of $\cM$,
is given by the formula
$$  \O_\cM \ = \ \I_{\cX/\cZ } \o_{FS}^{n+1}
\eqno(4.10)
$$
which is positive, since $\o_{FS}$ is positive.
Since $SU(N+1)$ is
semi-simple, there is a unique equivariant moment map
$$ \m_\cM: \cZ \ra su(N+1)
\eqno(4.11)
$$
with respect to the form $\O_\cM$ (see, for example, \S 4.9 of [CS]).

\v
\v
{\bf Theorem 1.}\ {\it Let $\cM$ be the line bundle over $\cZ$
given by {\rm (4.9)}. Then Donaldson's symplectic form
$\Omega_D$ and moment map $\m_D$ are given respectively by the
curvature $\Omega_{\cM}$ and the moment map $\m_{\cM}$ of $\cM$:
$$
\m_\cM=\m_\cD, \ \ \ \O_\cM = \O_D.
$$
}
\v
\v
{\it Proof.  } To prove this theorem, we compute $\m_\cM$ using
(3.4): Fix $z=[\us]\in \cZ$ and let $\g\in \cM$ be a point above $z$.
Let $\xi\in su(N+1)$ and define
$H(t)= -\log |\si_t\cdot\g|_\cM^2$ where
$\si_t=\exp(it\xi)$. Then the formula (3.4) can be rewritten as
$$
H'(0)\ = \ 2\langle \m_\cM(z), \xi\rangle
\eqno(4.12)
$$
On the other hand, the change of metric formula (2.7) for the Deligne
pairing tells us that
$$ H(t)-H(0) \ = \ E(\phi_t)\ = \
\I_X \phi_t\cdot \s_{j=0}^n (\si_t^*\o_{FS}^j\wedge \o_{FS}^{n-j} )
\eqno(4.13)
$$
where
$$ \phi_t(x) \ = \ \log{ |\si_t( x)|^2\over |x|^2 }
\eqno(4.14)
$$
The right hand side can be recognized as the familiar
component $F_{\o}^0(\phi_t)$ of the energy functional
$F_{\o}(\phi_t)$ due to Yau and Aubin.
Its variational derivative is well-known and
can be obtained by a straightforward computation
$$  {d\over dt}E(\phi_t)\ = \ \I_X \dot\phi_t \si^*_t\o_{FS}^n
\eqno(4.15)
$$
(see, for example, [S], Section \S 3.1).
Now $\dot\phi_t$ is given explicitly by
$$
\dot\phi_t\ = \ {2x^*\exp(2i\xi t) (i\xi)x\over  x^*\exp(2i\xi t)
x}
\eqno(4.16)
$$
for $x\in\P^N$. Differentiating (4.13) and substituting in (4.15) and
(4.16)
produces
$$
H'(0)\ = \ 2Tr\bigg(i\xi\I_X {xx^*\over x^*x} \o^n\bigg)\
\eqno(4.17)
$$
Comparing (4.12) and (4.17) we see that
$$ \m_\cM(z)\ = \ P_{su(N+1)}\left(i\I_X {xx^*\over x^*x} \o^n\right)\ = \
P_{su(N+1)}\left(i\langle s_\al,s_\b \rangle_{h_\us} \right)\ = \ \m_D(z)
\eqno(4.18)
$$
In view of the defining relation (3.1) between the moment map
and the symplectic form, and the fact that the symplectic forms
are compatible with the complex structure,
it follows that the
symplectic forms $\O_\cM$ and $\O_D$ must coincide on
each $SL(N+1)$ orbit. But then $\O_\cM=\O_D$
on $\cZ$, since
$\cZ$ consists of a single $SL(N+1)$ orbit.
\v\v

{\bf \S 5. Estimates for the moment map}
\v
In order to explain the statement of our theorem,
we sketch briefly the steps in Donaldson's proof [D],
introducing the necessary notation along the way.
One wants to show
 that constant scalar curvature
implies the existence of a ``balanced basis" of $H^0(L^k)$,
that is, an ordered basis $\us=\{s_0,..., s_N\}$ which imbeds $X$ into
$\P^N$ and which satisfies the following two conditions:

\vskip .03in
i) $\s_{j=0}^N |s_j|^2_{h_{FS}} \ = \ 1$ at each point of $X$.
\vskip .02in
ii) The matrix $\langle s_\al, s_\b\rangle_{h_{FS}}$ is diagonal,
where $h_{FS}$ is the Fubini-Study metric.

\v
To find the balanced basis, [D] first uses Lu's formula [L] for
the coefficients in the
Tian-Yau-Zelditch asymptotic expansion
for the density of states to produce, for each $k$, an ordered basis
$\us'=\{s_0',...,s_N'\}$ of $H^0(L^k)$ which is almost balanced in the
following sense:
$\s_{j=0}^N |s_j'|^2_{h_{FS}} \ = \ 1$ at each point of $X$ and
$$
\langle s_\al', s_\b'\rangle_{h_{FS}} = D_k + E_k
\eqno(5.1)
$$
where $D_k$ is a
scalar matrix with $D_k\ra 1$ as $k\ra\i$, and $E_k$ is a trace-free
hermitian
matrix whose operator norm $||E_k||_{op}$ tends to zero rapidly as $k$
tends
to infinity. He then interprets $E=E_k$ as the value $\m_D(z')$ where
$\m_D$ is the moment
map corresponding to  $\O_D$ on the manifold
$\cZ$, and $z'\in \cZ$ is the point
determined by $\us'$. The problem now becomes one of showing
that if $\m_D(z')$ is small for some $z'\in\cZ$, then there is a point
$z\in
\cZ$ which is close to
$z'$ with $\m_D(z)=0$. The standard technique for
finding the zero of a moment map is to follow the gradient flow
of the function $|\m_D|^2$. But to guarantee that the flow will,
after a short distance $\d$, reach a zero of the moment map, one
must prove that $d\m_D$,  the derivative of  $\m_D$,  is large in a
$\d$ neighborhood of $z'$.
\v
The heart of the argument in [D] is the proof of a
lower bound estimate on $d\m_D$. This is equivalent
to bounding $|\si_\cZ(\xi)|$ from below, where $\xi\in su(N+1)$
has length one and $\si_\cZ(\xi)$ is the vector field on $\cZ$
given by the infinitesimal action of $\xi$.
Since $\cZ\sub\cH_0//\cG$,
the  symplectic reduction formalism tells us that $\si_\cZ(\xi)$ is
the projection of $\si_{\cH_0}(\xi)$ (the vector field on
$\cH_0$ corresponding to the infinitesimal action of $\xi$) onto the
orthogonal complement of the tangent space $T(\cG^c\cdot w)$, where
$\cG^c\cdot w$ is a complexified orbit. This description turns out to
be sufficiently explicit to allow the following estimation of
$|\si_\cZ(\xi)|$:
\v
\v
{\bf Theorem.} (Donaldson) {\it\
Suppose $Aut(X,L)$ is discrete. For any $R>1$ there
are positive constants  $C$ and $\e<{1\over 10}$ such that,
for any $k$, if the basis $s_\al$ of $H^0(\cL^k)$
has R-bounded geometry, and if $||E||_{op}< \e$, then
$$
\Lambda_z \ \leq \ C^2\cdot k^4
\eqno(5.2)
$$
}
\v
Here the basis $s_{\al}$ is viewed as defining a metric
$\tilde\o$ on $X$, which is in the cohomology class $k\,c_1(L)$.
The metric $\tilde\o$
is said to have $R$-bounded geometry if $\tilde\o>R^{-1}\tilde\o_0$
and $||\tilde\o-\tilde\o_0||_{C^r(\tilde\o_0)}<R$,
where $\tilde\o_0=k\o_0$, $\o_0$ is a fixed reference metric in $c_1(L)$,
$||\cdot||_{C^r(\tilde\o_0)}$ is the $C^r$ norm defined by
$\tilde\o_0$, and $r\geq 4$ is a fixed integer.
The basis $s_{\al}$ is said to have $R$-bounded geometry
if the corresponding metric $\tilde\o $ does.
The expression
$\Lambda_z^{-1}$ denotes the smallest eigenvalue of the
operator
$$
Q_z=\sigma_z^*\sigma_z:\ su(N+1)\ \longrightarrow\ su(N+1),
$$
where $\sigma: su(N+1)\ra T\cZ$
is the infinitesimal action and $\sigma_z^*$ is its
adjoint with respect to the metric on $T\cZ$ and the
invariant metric on $su(N+1)$. At the end of his paper,
Donaldson sketches a
refinement of the argument which replaces $k^4$ by $k^{2+\e}$, and
conjectures that in fact the estimate (5.2) holds
with $k^4$ replaced by $k$.
In this section, we shall prove:

\v\v

{\bf Theorem 2.} {\it Under the same hypotheses as in Donaldson's theorem,
we have}
$$
\Lambda_z\ \leq \ C^2\cdot k^2
\eqno(5.3)
$$
\v\v
{Proof. } The  estimate (5.3) is equivalent to the estimate
$$
|\sigma_{\cZ}(\xi)|^2\ \geq\ c_Rk^{-2}||\xi||^2
\eqno(5.4)
$$
for all $\xi\in su(N+1)$
and a positive constant $c_R$ depending only
on $R$. To prove this, our starting point is the identity
$$
|\si_\cZ(\xi)|^2 \ = \ \I_X \iota_{X_\xi,\bar X_\xi }\o_{FS}^{n+1}
\eqno(5.5)
$$
which follows from Theorem 1. Indeed,
$|\si_\cZ(\xi)|^2=-i\Omega_D(\si_\cZ(\xi),\overline{\si_\cZ(\xi)})$
by definition. On the other hand,
$\Omega_D(\si_\cZ(\xi),\overline{\si_\cZ(\xi)})
=
\Omega_\cM(\si_\cZ(\xi),\overline{\si_\cZ(\xi)})$
by Theorem 1, and
this last expression is given
by the right hand side of (5.5) in view
of Deligne's formula (2.9).
(Alternatively, (5.5) has also been proved directly in
[PS], formula (3.5)).
\v
Consider now the exact sequence of holomorphic vector bundles:
$$
0 \ra TX \ra \iota^*T\P^N \ra Q \ra 0
$$
where $TX$ is the tangent bundle of $X$, $T\P^N$ is  the
tangent bundle of $\P^N$, $\iota:X\ra\P^N$ is the
embedding, and $Q$ is the quotient.
Let
$\cN\sub \iota^*T\P^N$ be the orthogonal complement of
$TX$.  Then $\cN$ is a
smooth vector bundle and
$$
\iota^*T\P^N\ = \ TX\oplus \cN
$$
Let $\pi_T$ and $\pi_\cN$ be the projections onto
the first and second components of $\iota^*T\P^N$. We observe that
$$
\I_X \iota_{X_\xi,\bar X_\xi }\o_{FS}^{n+1}
=
||\pi_\cN X_{\xi}||^2
\eqno(5.6)
$$
where $||\cdot||$ denotes the $L^2$ norm with respect to
the metric $\tilde\o$ on the base and the Fubini-Study metric
on the fiber.
The desired inequality (5.4) is a direct consequence of (5.5), (5.6)
and the following inequalities:
$$
\eqalignno{
||\xi||^2&\leq c_R'
   k||X_{\xi}||^2 \
&(5.7)\cr
||X_\xi||^2\ &= \ ||\pi_TX_\xi||^2+||\pi_NX_\xi||^2\
&(5.8)\cr
c_R\,||\pi_TX_{\xi}||^2&\leq \ k||\pi_\cN X_\xi||^2
&(5.9)\cr}
$$
where $c_R,c_R'>0$  are constants independent of $k$.
\v
The identity (5.8) is an immediate consequence of the orthogonality of
$\pi_NX_\xi $ and $\pi_TX_\xi $.
To prove (5.7) we argue as follows:
If we view points
$z\in \P^N$ as a column vectors in $\C^{N+1}$ (modulo the action of
$\C^\times$), then tangent vectors
$X$ on
${\bf CP}^N$ can be viewed as ordered pairs $(z,v)$ modulo the equivalence
relation:
$(z,x)\sim (z',v')$ if $z'=\lambda z$, $v'-\lambda v=\mu z$ for
some $\l\in\C^\times$ and $\m\in\C$.
Then the vector field $X_{\xi}$ is
defined, for $\xi\in su(N+1)$,  by $X_{\xi}=(z,\xi z)$, and its
norm with respect to the Fubini-Study metric $\omega_{FS}$ is given by the
well-known formula
$$
|X_\xi|^2(z)\ = \ {(z^*\xi^*\xi z)(z^*z)-(z^*\xi z)^2\over (z^*z)^2}
$$
Observe that the function $z^*\xi z/z^*z$ on $M$ is just
$\dot\phi\equiv\dot\phi_t$, at $t=0$, where $\phi_t= \log{|\si_t z|^2\over
|z|^2}$ and
$\si_t(z) = \exp(it\xi)z$.
Integrating the above identity with respect to $\omega_{FS}$ gives
$$
tr\left(
\xi^*\xi\cdot
\I_M {zz^*\over z^*z}\tilde\o^n\right)
=||X_{\xi}||^2+\int_M{\dot\phi_t^2}\tilde\o^n
\eqno(5.10)
$$

We claim that the following
Poincare inequality holds for $\tilde\omega$,
with uniform constants in $\tilde\o$
$$
c\int_M\dot\phi^2\tilde\omega^n
\leq k\int_M\bar\pl\dot\phi\wedge\pl\dot\phi\wedge\tilde\omega^n
+k^{-n}(\int_M\dot\phi\,\tilde\omega^n)^2
\eqno(5.11)
$$
Indeed, for $\tilde\omega_0\equiv k\o_0$, this is just the
standard Poincar\'e inequality for the fixed metric $\o_0$, scaled up by a
factor of $k$. To establish it for $\tilde\o$, we note first that
$R^{-1}\tilde\omega_0<\tilde\omega<2R\tilde\omega$ since $\tilde\omega$ is
$R$-bounded.
Writing
$$
\int_M\dot\phi\,\tilde\omega_0^n
=
\int_M\dot\phi\,\tilde\omega^n
-
\int_M\dot\phi(\tilde\o^n-\tilde\o_0^n)
$$
and $\tilde\o^n-\tilde\o_0^n=\pl\bar\pl\,\theta\wedge
\sum_{p=0}^{n-1}\tilde\o_0^{n-1-p}\tilde\o^p$, with $\t$ normalized
so that $\I_M\t\tilde\o_0^n=0$, we have
$$
\eqalign{|\int_M\dot\phi(\tilde\o^n-\tilde\o_0^n)|
&=
|\int_M\bar\pl\dot\phi\wedge\pl\theta\wedge
\sum_{p=0}^{n-1}\tilde\o_0^{n-1-p}\tilde\o^p|
\leq
\sum_{p=0}^{n-1}
\int_M|\bar\pl\dot\phi|_{\tilde\o_0}|\pl\theta|_{\tilde\o_0}
({\tilde\o\over\tilde\o_0})^p\tilde\o_0^n\cr
&\leq
C_1(\int_M|\bar\pl\dot\phi|_{\tilde\o_0}^2\tilde\o_0^n)^{1/2}
(\int_X|\bar\pl\dot\theta|_{\tilde\o_0}^2\tilde\o_0^n)^{1/2}\cr
&\leq C_2(\int_M\bar\pl\dot\phi\wedge\pl\dot\phi\wedge
\tilde\o_0^{n-1})^{1/2}k^{{1\over 2}(n+1)}
\cr}
$$

The last inequality follows from the
uniform boundedness of $||\D_{\tilde\o_0}\t||_{C^\i}$
(which holds by the R-boundedness assumption)
and the inequality
$||\bar\pl\t||^2_{\tilde \o_0} \leq Ck||\D_{\tilde\o_0}\t||^2_{\tilde
\o_0}$ (which holds since $k\Delta_{\tilde \o_0}-c>0$
on the space of $\theta$'s
with mean value $0$,
for some small positive constant $c$). The inequality (5.11)
now follows.
\v
Writing
$\I_M {zz^*\over z^*z}\tilde\o^n\ = \ D_k+E_k
$
and using the fact that $D_k\ra 1$ as $k\ra \i$, and $||E_k||_{op}<\e$,
we see that
$
tr\left(
\xi^*\xi\cdot
\I_M {zz^*\over z^*z}\tilde\o^n\right)\geq c||\xi||^2.
$
On the other hand,
$$|\int_M{\dot\phi_t}\tilde\o^n|=|tr(\xi E)|\leq
\sqrt{N+1}||\xi||\cdot  ||E||_{op}\ \leq ck^{n/2}||\xi||\cdot ||E_k||_{op}
$$
Combining this with (5.10) and (5.11), and using $||E_k||_{op}<\e$, we
obtain, for $\e$ small,
$$ c||\xi||^2\ \leq \ ||X_\xi||^2 +
k\int_M\bar\pl\dot\phi\wedge\pl\dot\phi\wedge\tilde\omega^n
$$
But now observe that $\bar\pl\dot\phi=\iota_{X_\xi}\o_{FS}$. Restricting
to $M$ we get $\bar\pl\dot\phi=\iota_{\pi_T{X_\xi}}\tilde\o$ which implies

$$ c||\xi||^2\ \leq \ ||X_\xi||^2 +
k||\pi_TX_\xi||^2
$$
and this proves (5.7).

\v
It remains to prove
(5.9).
We begin by assuming that $\Aut(X)$ is discrete, that is, $X$
does not admit any nonzero holomorphic vector field. This assumption
will be removed later.
Since there are no nonzero holomorphic vector fields on $X$,
we have the
standard inequality for the $\bar\partial$ operator:
$$
c\,||W||_{L^2(\o_0)}^2 \ \leq \ \,||\bar\pl(W)||_{L^2(\o_0)}^2
\eqno(5.12)
$$
for some $c>0$,
where $W$ is any smooth vector field on $X$ and $||\cdot ||_{L^2(\o_0)}$
is the $L^2$ norm with respect to the metric $\o_0$.
Replacing $\o_0$ by $\tilde\o_0=k\o_0$, we obtain the inequality
$ c\,||W||_{L^2(\tilde\o_0)}^2 \ \leq \ \, k\,
||\bar\pl(W)||_{L^2(\tilde\o_0)}^2
$\ .
Recall that $||\cdot||$ denotes $L^2$-norms with respect to the metric
$\tilde\o$ on the base and the Fubini-Study metric on
the fiber. Since $\tilde\o$ has $R$-bounded geometry, we have
$$ c_R||W||^2 \ \leq \ \, k\,
||\bar\pl(W)||^2
\eqno(5.13)
$$
for some $c_R>0$, depending on $R$ but not on $k$.
Applying this inequality with $W=\pi_T X_\xi$,
we see that in order to establish (5.11) it suffices to prove
$$
c_R\,||\bar\pl(\pi_TV)||^2\
\leq
||\pi_\cN V||^2,
\eqno(5.14)
$$
for all holomorphic vector fields $V$ on $\P^N$.
\v

But $V=\pi_T V + \pi_\cN V$ so $\bar\pl V = 0  =
\bar\pl(\pi_TV)+\bar\pl(\pi_\cN V)$. Thus
it suffices to prove:
$$ ||\pi_\cN V||^2 \ \geq \ {c_R}\cdot ||\bar\pl(\pi_\cN V)||^2
\ = \ {c_R}\cdot ||(\bar\pl\pi_\cN)( V)||^2
\eqno(5.15)
$$
In fact, we shall prove the pointwise estimate:
$$ |\pi_\cN V|^2 \ \geq \ {c_R}\cdot |\bar\pl(\pi_\cN V)|^2
\ = \ {c_R}\cdot |(\bar\pl\pi_\cN)( V)|^2
\eqno(5.16)
$$
Fix $x\in X$ be a point and choose a local holomorphic frame of
$\iota^*T\P^N$, in a neighborhood of $x$, of the form: $e_1,...,e_n,
f_1,...,f_m$ with
$m+n=N$,  satisfying the following:
\vskip .03in
a) The $e_1,...,e_n, f_1,...,f_m$ form an orthonormal basis
of $\iota^*T_x\P^N$.
\vskip .03in
b) The $e_1,...,e_n$ form a local holomorphic frame of $TX$ near $x$.
\v
Then we can express $V=\sum_ia_ie_i+\sum_j b_jf_j$,
where the $a_i$ and the $b_j$ are holomorphic functions.
Since $\pi_\cN$ is a projection, $\pi_\cN(\pi_\cN(f_j)-f_j)=0$.
Thus $\pi_\cN(f_j)-f_j$ is a linear combination of the vectors
$e_i$'s, and we can write
$$ \pi_\cN(f_j)\ = \ f_j - \sum_{i=1}^n \phi_{ij}e_i
\eqno(5.17)
$$
where the $\phi_{ij}$'s are smooth functions, vanishing at $x$. Then
$\pi_\cN V=\s_{j=1}^m b_j \big(f_j-\s_i \phi_{ij}e_i\big)$, and
$$
\bar\pl(\pi_\cN V) \ = \ \s_{j=1}^m b_j \big(-\s_{i=1}^n
(\bar\pl\phi_{ij})e_i\big)
\eqno(5.18)
$$
To establish (5.16) we must prove:
$$ \s_{i=1}^n |\s_{j=1}^m b_j\bar\pl\phi_{ij}|^2\ \leq
c_R^{-1}\cdot \s_{j=1}^m|b_j|^2
\eqno(5.19)
$$
But for each $i$, the Cauchy-Schwarz inequality
for quadratic forms defined by hermitian matrices implies
$|\s_{j=1}^m b_j\bar\pl\phi_{ij}|^2\ \leq \ \sum_{j=1}^m |b_j|^2\cdot
\sum_{j=1}^m |\bar\pl\phi_{ij}|^2$,
and thus it suffices to prove
$$
\sum_{i=1}^n\sum_{j=1}^m |\bar\pl\phi_{ij}|^2\ \leq \ c_2
\eqno(5.20)
$$
where $c_2=c_2(R)$ is  independent of $k$ (and depends only on $R$).
But the matrix
$$
A^* \ =\ (\bar\pl\phi_{ij})
\eqno(5.21)
$$
is precisely the dual of the second fundamental form $A$
of $TX$ in
$\iota^*T\P^N$. One can see this as follows: Let $S=TX$,
$E=\iota^*T\P^N$ and $Q=E/S$. Then we have the exact sequence
of holomorphic bundles:
$$  0 \ra S \ra E \mapright p Q \ra 0
\eqno(5.22)
$$
The second fundamental form is the bundle map
$A: S \ra Q\otimes \O^{1,0}$ defined by $A=p\circ D_E$, where
$D_E$ is the Chern connection on $E$, compatible with the metric
$h=h_{FS}$ and
the holomorphic structure. Then $A^*: Q\ra S\otimes \O^{0,1}$
is characterized by
$$ \langle As,q\rangle_h\ = \ \langle s, A^*\tilde q\rangle_h\ = \
\langle s, A^* q\rangle_h\
\eqno(5.23)
$$
where $s$ is a section of $S$ and $q$ a section of $Q$, and
$\tilde q$ is the section of $S^\perp$ which corresponds to $q$
via the canonical isomorphism of smooth bundles $S^\perp\ra Q$.
We claim that
$$
A^*(\tilde f_j)=\s_{i=1}^n(\db \phi_{ij})e_i
\eqno(5.24)
$$
To see this, first note that $\langle s, \tilde f_j\rangle= $
$\langle s, f_j-\s_{i=1}^n\phi_{ij}e_i\rangle=0$.
Differentiating both sides:
$$  \langle D_E s, f_j-\s_{i=1}^n\phi_{ij}e_i\rangle +
\langle s, D_E(f_j-\s_{i=1}^n\phi_{ij}e_i)\rangle = 0
\eqno(5.25)
$$
In other words,
$$  \langle As, f_j-\s_{i=1}^n\phi_{ij}e_i\rangle\ \ = \
\langle s, \s_{i=1}^n \db(\phi_{ij})e_i\rangle
\eqno(5.26)
$$
This proves the claim. Since $A$ is the second fundamental
form, the curvature tensors  of the sub-bundle and of the ambient
bundle are related by (see, e.g.,  [GH], page 78):
$$
-A^*\wedge A\ = \ \pi_T\circ (F_{\iota^*T\P^N}|_{TX})\ -\ F_{TX}
\eqno(5.27)
$$
Let us
explain the meaning of this notation: The tensor
$F_{\iota^*T\P^N}$ is a 2-form on $X$ with values in $End(\iota^*T\P^N)$.
Thus
$F_{\iota^*T\P^N}|_{TX}$ is a  2-form with values in
$Hom(TX, \iota^*T\P^N)$ and
$\pi_T\circ (F_{\iota^*T\P^N}|_{TX})$ is a 2-form with values in
$Hom(TX,TX)=End(TX)$ as is the tensor $F_{TX}$.
Thus
$$
\sum_{i=1}^n\sum_{j=1}^m |\bar\pl\phi_{ij}|^2\ = \
C_{\tilde\o}Tr\big[(\pi_T\circ (F_{\iota^*T\P^N}|_{TX})-F_{TX})
\big]
\eqno(5.28)
$$
where
$$
Tr\big[(\pi_T\circ (F_{\iota^*T\P^N}|_{TX})-F_{TX})
\eqno(5.29)
$$
is a two form on $X$ and
$C_{\tilde\o}$ is the contraction with the metric $\tilde\o$.
\v
We now compute the terms in the equation above: First,
$R=F_{T\P^N}$ is the Riemann curvature of the Fubini-Study
metric. It is well known that the Fubini-Study metric
has constant bisectional curvature. In fact:
$$
R(X,\bar Y, Z, \bar W) \ = \ g(X,\bar Y)g(Z,\bar W)+ g(X,\bar W)
g(Z,\bar Y)
\eqno(5.30)
$$
where $g=g_{FS}$ is the Fubini-Study metric. In otherwords,
in any system of local coordinate for $\P^N$, we have
$R_{i\bar j k \bar l} \ = \ g_{i\bar j}g_{k\bar l}+g_{i\bar l}
g_{k\bar j}$, and thus
$$
{{R_i}^j}_{k\bar l} \ = \ \d_i{}^jg_{k\bar l}+ g_{i\bar l}\d^j{}_k
\eqno(5.31)
$$
Now for a fixed $x\in X$, choose a local coordinate system
$(x_1,.., x_N)$, centered at $x$,
in such a way that $T_xX$ is the plane defined by:
$x_{n+1}=\cdots = x_N = 0$. Then
$$ Tr\big[(\pi_T\circ (F_{\iota^*T\P^N}|_{TX})\big]\ = \
\s_{k,l=1}^n\left(\sum_{i=1}^n {{R_i}^i}_{k\bar l}\right)
dz_k\wedge d\bar z_{ l}
\  \
$$
$$ =\ n \s_{k,l=1}^ng_{k\bar l}\cdot dz_k\wedge d\bar z_{ l} \ + \
\s_{k, l=1}^n g_{k\bar l}\cdot dz_k\wedge d\bar z_{l}\ = \
(n+1)\tilde\o
\eqno(5.32)
$$
\v
Thus
$$
C_{\tilde\o}Tr\big[(\pi_T\circ (F_{\iota^*T\P^N}|_{TX})\big]\ = \ (n+1)
\eqno(5.33)
$$
\v

On the other hand,
$C_{\tilde\o}(F_{TX}) $ is the scalar curvature of $X$ with
respect to the pullback of $\o_{FS}$. But our assumption in
Theorem 2 is that $\o_{FS}$ has $R$-bounded geometry with
respect to $k\o_0$. This is readily seen to imply that
$||\nabla^r\tilde\o||_{C^0(\o_0)}\leq C_R\,k^{1+{r\over 2}}$,
and hence
$$
C_{\tilde\o}Tr\big[F_{TX}\big]\ \leq \ C_R
\eqno(5.34)
$$
with a constant depending
only on $R$. This completes the proof under the assumption
$Aut(X)$ is discrete.
\v
Finally, we remove the hypothesis that $Aut(X)$ is discrete.
This step relies on Lemma 12 of [D], which we now review:
\v
Let $p: L\ra X$ be a positive holomorphic vector bundle on a compact
manifold
$X$. Let $Aut(X)$ be the group of holomorphic automorphisms of $X$ and
let $Aut(X,L)$ the  group of holomorphic automorphisms of the
pair $(X,L)$. Thus an element of $Aut(X,L)$ is a pair $(F,\hat F)$
where $F:X\ra X$ is biholomorphic, $\hat F: L\ra L$ is biholomorphic,
and the diagram commutes: $p\hat F = Fp$.
\v
We clearly have a map $Aut(X,L)\ra Aut(X)$ defined by
$(F,\hat F)\mapsto F$. The kernel of the map is $\C^\times$.
We are interested in  the image of this homomorphism.
More precisely, we want to characterize the infinitesimal image, that is,
the image of $Lie(Aut(X,L))$ inside $Lie(Aut(X))$.
\v
Recall that
if $v$ is a vector field on $X$ then
$v\in Lie(Aut(X))$ if and only if $v$ is the real part
of a holomorphic vector field, i.e.:
\vskip .03in
1. $v= w+\bar w$ where $w=v^j{\pl\over dz^j}$

2. $\bar\pl w=0 $
\v
 An element of
$Lie(Aut(X,L))$ is a vector field $V$ in $Aut(L)$
(i.e., $V$ is the real part of a holomorphic vector
field on $L$) which
is $\C^\times$ invariant. We have a well defined
map  $q: Lie(Aut(X,L))\ra Lie(Aut(X))$: If $V\in Lie(Aut(X,L))$,
and if $x\in X$, then $q(V)(x)\in T_xX$ is defined
as follows: Let $l\in L$ be any point in $p^{-1}(x)$.
Then $q(V)(x)= dp(V(l))$. This does not depend on
the choice of $l$, since $V$ is $\C^\times$ invariant.
\v

The kernel
of $q$ is $\C\cdot\bt$, where $\bt$ is the vector field
generated by the infinitesimal action of $U(1)$. The hypothesis
which Donaldson imposes in his theorem is that
the image of $q$ is trivial, that is,
$Lie(Aut(X,L))/(\C\cdot \bt)$ is trivial.  This group
is characterized in Lemma 12 of [D] as follows:
Fix a hermitian
metric $h$ on $L$ with positive curvature $\o$ (such an $h$
exists since, by assumption, $L$ is positive). Then Lemma 12 of [D] says
that  $v$ is in the
image of $q$ if and only if $w^i= \o^{i\bar j} f_{\bar j}$ for some
smooth
complex valued function $f$ on $X$.
Note that this does not depend on the choice of $h$: If $h'$ is
another positive metric, then $\o$ is replaced by
$\o'=\o+{\sqrt{-1}\over 2\pi}\ddb \phi$ for some $\phi$. So if
$\o_{i\bar j}w^i = \db f$
then $\o'_{i\bar j}w^i =
\db f + {\sqrt{-1}\over 2\pi}\phi_{i\bar j}w^i = \db f'$
where $f'=f+ {\sqrt{-1}\over 2\pi}\phi_iw^i$ (here we are using the fact
that $\db w=0$).
Note as well that if we replace $L$ by $L^k$, the image of $q$ does
not change since if $\o_{i\bar j}w^i = \db f$ then
$k\o_{i\bar j}w^i = \db (kf)$.
\v
We now show how Lemma 12 of [D] can be used to remove the
hypothesis that $Aut(X)$ is discrete:
Let $L\ra X$ be a hermitian line bundle on a compact
complex manifold $X$. Assume that $\o$, the curvature
of $L$, is a positive (1,1) form.
If $X$ has
no holomorphic vector fields, then there is a
constant $c>0$ such that
$$  || w ||_{L^2_1(\o_0)} \ \leq \ c\cdot || \db w||_{L^2(\o_0)}
\eqno(5.35)
$$
for all $L^2_1$ vector fields $w=w^i{\pl\over \pl z_i }$ on $X$.

\v
Now we drop the assumption that $X$ has no holomorphic
vector fields. Then (5.35) no longer holds for all
$w\in L^2_1$. Let $H_0\sub L^2_1$ be the space of vector fields
of the form $w^i = \o^{i\bar j}f_{\bar j}$ where
$f$ ranges over all smooth complex valued functions
on $X$, and let $H\sub L^2_1$ be the completion of $H_0$.
By the standard elliptic estimates for the $\bar\partial$
operator, the space $H$ can be identified with vector fields
of the form $\o^{i\bar j}f_{\bar j}$, where $f$ is a function
in the Sobolev space $L_2^2$.
\v
{\bf Lemma : } Assume that $Aut(X,L)/\C^\times$ is discrete. Then
(5.35) holds for all $w\in H$.
\v
Proof. Lemma 12 of [D] says that the kernel of $\db|_{H_0}$
  is
$Lie(Aut(X,L))/\C$. But by elliptic regularity, $ker(\db|_H)=
ker(\db|_{H_0})$.  Now we are assuming
$$
Lie(Aut(X,L)/\C^\times)\ = \ Lie(Aut(X,L))/\C\ = \ 0
\eqno(5.36)
$$
Thus the kernel of $\db|_H$ is
trivial. Now let $w\in H$. Then $w'=w-\t\in ker(\db)^\perp$,
for
some unique $\t\in ker(\db)$. Here the orthogonal complement
$ker(\db)^{\perp}$ is taken with respect to
the $L^2$ norm. Thus
$$  || w' ||_{L^2_1(\o_0)} \ \leq \ c\cdot || \db w'||_{L^2(\o_0)}\ = \
c\cdot || \db w||_{L^2(\o_0)}
\eqno(5.37)
$$
On the other hand, there is a constant $c'>0$ such that
$$
||w||_{L^2(\o_0)} \ \leq \ c'||w'||_{L^2_1(\o_0)}
\eqno(5.38)
$$
for all $w\in H$. To see this, assume that $||w'_n||_{L^2_1(\o_0)}\ra 0$
for some sequence $w_n\in H$ such that $||w_n||_{L^2(\o_0)}=1$.
 Since
$w_n=w_n'+\t_n$ is an orthogonal decomposition in $L^2$,
we have $||\t_n||_{L^2(\o_0)}\leq 1$. But the $\t_n$ are holomorphic
vector fields. Thus, after passing to a subsequence, the
$\t_n$ converge in $C^\i$ to an element $\t\in ker(\db)$. Since
$w_n'=w_n-\t_n\ra 0$ in $L^2_1(\o_0)$, we conclude that
$w_n \ra \t$ in $L^2_1(\o_0)$.  Since $w_n\in H$ and since
$H$ is closed in $L^2_1(\o_0)$, we have $\t\in H$. But we have
seen that $ker(\db|_H)=0$. Thus $\t=0$. On the other hand,
$w_n\ra \t$ in $L^2_1(\o_0)$ so $1=||w_n||_{L^2(\o_0)}\ra 0$, which is
a contradiction. This proves (5.38), and together with (5.37),
$||w||_{L^2(\o_0)}\leq c'||\db w||_{L^2(\o_0)}$. Our claim follows,
in view of the standard
elliptic estimate
$$
||w||_{L_1^2(\o_0)}\leq c''(||\db w||_{L^2(\o_0)}+||w||_{L^2(\o_0)}).
\eqno(5.39)
$$
Now we conclude the proof of Theorem 2.
The discreteness assumption of $Aut(X)$ was used
 earlier only at one point, namely  to ensure that (5.12)  holds
with $W=\pi_TV$.
The  lemma shows that if we only assume that $Aut(X,L)/\C^\times$ is
discrete, that
this inequality continues to hold provided $W$ is in $H$.
Thus it suffices to show that
$\pi_TV\in H$:
Let $V$ be the vector field on $\P^N$ given by the infinitesimal
action of some element in $sl(N+1)$.
Then $\o^{FS}_{i\bar j} V^i=\db F$ for some $F$, smooth on $\P^N$.
In other words,
$ \db F(Y)\ = \ \o^{FS}(V,Y)
$
for any $Y\in T\P^N$. Now suppose that $Y\in TX$, and let
$f=F|_X$ and let $\o = \o^{FS}|_X$.
Then
$$
\db_X f(Y) = \db_{\P^N} F(Y)\ = \ \o^{FS}(V,Y)\ = \ \o^{FS}(\pi_T V,Y)\
=
\
\o(\pi_T V,Y)
\eqno(5.40)
$$
so
$$ (\pi_T V)^i\ = \ \o^{i\bar j}f_{\bar j}
\eqno(5.41)
$$
Thus $\pi_T V\in H_0\sub H$, which is what we wanted to prove.

\v
{\it Remark 1.} Donaldson also conjectures that
$\Lambda_z\leq Ck
$
 holds with the
Hilbert-Schmidt norm on $su(N+1)$ replaced by the operator norm. This
problem
is still open.
\v
{\it Remark 2.} Let $X=\P^1$ be imbedded in $\P^k$ by $O(k)$. Consider
$\xi\in su(k+1)$ given by
$$ \xi_{ii}\ = \ i^2-{1\over 6}k(2k+1)-k(i-{1\over 2}k)
$$
and $\xi_{ij}=0$ if $i\not= j$. Then $\xi$ is orthogonal to the Lie
algebra of $ \Aut (\P^1,O(k))$ and
$$ ||\xi||^2\ = \ {1\over 180}k^5 + O(k^4),\hskip .2in||X_\xi||^2\ = \
{1\over
30}k^4 + O(k^3)
\hskip .2in||\pi_TX_\xi||^2\ = \ {1\over
30}k^4 + O(k^3)
$$
Thus $||\pi_NX_\xi||^2=O(k^3)$, and the bound $||\xi||^2\leq
C\,k^2||\pi_NX_\xi||^2 $ is sharp.
\v\v
{\bf Acknowledgements}: We would like to thank the referees and Xiaowei
Wang for pointing out to us a gap in our preprint.
\vfill\break

\centerline{\bf REFERENCES}

\v\v

[CS] Cannas da Silva, A.,
{\it Lectures on symplectic geometry},
Lecture Notes in Math. {\bf 1764} (2001)
Springer-Verlag, Berlin.

\v

[C] Catlin, D., ``The Bergman kernel and a theorem of Tian",
{\it Analysis and geometry in several complex variables,
Katata, 1999}, eds. G. Komatsu
and M. Kuranishi,
Trends in Math. (1999) 1-23,
Birkh\"auser, Boston.
\v

[De] Deligne, P.,
``Le determinant de la cohomologie",
Contemporary Math. {\bf 67} (1987) 93-177.

\v
[D] Donaldson, S.,
``Scalar curvature and projective imbeddings I",
J. Differential Geometry {\bf 59} (2001) 479-522.

\v
[DK] Donaldson, S. and P. Kronheimer,
{\it The Geometry of four-manifolds},
Oxford University Press, 1990.

\v
[GH] Griffiths, P. and J. Harris,
{\it Principles of Algebraic Geometry},
Wiley, 1978.

\v

[Lu] Lu, Z.,
``On the lower order terms of the asymptotic expansion
of Tian-Yau-Zelditch", Amer. J. Math. {\bf 122} (2000)
235-273.

\v
[L] Luo, H.,
``Geometric criterion for Gieseker-Mumford stability
of polarized manifolds",
J. Diff. Geom. {\bf 49} (1998) 577-599.

\v

[PS] Phong, D.H. and J. Sturm,
``Stability, energy functionals, and
K\"ahler-Einstein metrics",
Comm. in Analysis and Geom. {\bf 11} (2003)
563-597, arXiv:math.DG/0203254.

\v
[S] Siu, Y.T.,
``{\it Lectures on Hermitian-Einstein metrics for stable
bundles and K\"ahler-Einstein metrics},
Birkh\"auser, Boston, 1987.
\v

[T1] Tian, G.,
``On a set of polarized K\"ahler metrics
on algebraic manifolds",
J. Diff. Geometry {\bf 32} (1990)
99-130.
\v
[T2] Tian, G.,
``The K-energy on hypersurfaces and stability",
Comm. Anal. Geometry {\bf 2} (1994) 239-265.

\v

[T3] Tian, G.,
``K\"ahler-Einstein metrics with positive scalar
curvature", Inventiones Math. {\bf 130} (1997) 1-37.

\v
[T4] Tian, G.,
``Bott-Chern forms and geometric stability",
Discrete Contin. Dynam. Systems {\bf 6} (2000)
211--220.

\v
[Y1] Yau, S.T.,
``On the Ricci curvature of a compact K\"ahler manifold
and the complex Monge-Ampere equation I",
Comm. Pure Appl. Math. {\bf 31} (1978) 339-411.

\v

[Y2] Yau, S.T.,
``Open Problems in Geometry",
Proc. Symposia Pure Math. {\bf 54} (1993) 1-28.

\v
[Y3] Yau, S.T.,
``Review of K\"ahler-Einstein metrics in Algebraic Geometry",
Israel Math. Conf. Proceedings {\bf 9} (1996)
433-443.

\v
[Y4] Yau, S.T.,
``Nonlinear analysis in geometry",
Enseign. Math. {\bf 33} (1987) 109-158.
\v

[Ze] Zelditch, S.,
``Szeg\"o kernel and a theorem of Tian",
Int. Math. Research Notices {\bf 6} (1998) 317-331.
\v

[Z] Zhang, S.,
``Heights and reductions of semi-stable varieties",
Compositio Math. {\bf 104} (1996) 77-105.

\end